\documentclass[12pt]{article}
\usepackage{}
\usepackage{mathrsfs}

\usepackage{amsmath, epsfig, cite}
\usepackage{amssymb}
\usepackage{amsfonts}
\usepackage{latexsym}
\usepackage{graphicx}

\newtheorem{thm}{Theorem}[section]

\newtheorem{cor}[thm]{Corollary}
\newtheorem{defi}[thm]{Definition}
\newtheorem{lem}[thm]{Lemma}

\newcommand{\qed}{{\hfill\rule{4pt}{7pt}}}
\def\pf{\noindent {\it Proof.} }

\numberwithin{equation}{section}

\makeatletter \@addtoreset{equation}{section} \makeatother

\begin{document}
\rule{0cm}{1cm}
\begin{center}
{\Large\bf The matching energy of graphs with given edge connectivity\footnote{The first author is supported by  NNSFC (Nos. 11326216 and 11301306);
 the second author is supported by  NNSFC (Nos. 11101351 and 11171288) and
 NSF of the Jiangsu Higher Education Institutions (No. 11KJB110014)
 .}
}
\end{center}
 \vskip 2mm \centerline{ Shengjin Ji$^{1}$, Hongping Ma$^{2}$
 \footnote{Corresponding author.\\
E-mail addresses: jishengjin2013@163.com(S.Ji), hpma@163.com(H.Ma)}}

\begin{center}
$^{1}$ School of Science, Shangdong University of Technology,
\\ Zibo, Shandong 255049, China\\
$^{2}$ School of Mathematics and Statistics, Jiangsu Normal
University,\\
Xuzhou, Jiangsu 221116, China

\end{center}

\begin{center}
{\bf Abstract}
\end{center}

{\small Let G be a simple graph of order $n$ and
$\mu_1,\mu_2,\ldots,\mu_n$ the roots of its matching polynomial. The
matching energy of $G$ is defined as the sum $\sum_{i=1}^n|\mu_i|$.
Let $K_{n-1,1}^k$ be the graph obtained from $K_1\cup K_{n-1}$ by
adding $k$ edges between $V(K_1)$ and $V(K_{n-1})$. In this paper,
we show that $K_{n-1,1}^k$ has maximum matching energy among all
connected graph with order $n$ and edge connectivity $k$.}

\vskip 3mm

 \vskip 3mm


\section{Introduction}

We use Bondy and Murty \cite{graphBondy2008} for terminology and
notations not defined in this paper and consider undirected and
simple graphs only. Let $G=(V, E)$ be such a graph with order $n$.
 Denote by $m(G,t)$ the number of $t$-matchings of $G$. Clearly,
 $m(G,1)=e(G)$, the size of $G$, and $m(G,t)=0$ for $t > \lfloor n/2\rfloor$. It is both
 consistent and convenient to define $m(G,0)=1$.

 Recall that the \emph{matching polynomial} of a graph $G$ is defined as
$$\alpha(G)=\alpha(G,\lambda)=\sum\limits_{t\geq 0}(-1)^t m(G,t)\lambda^{2t}$$
and its theory is well elaborated \cite{CDGT,Farrell1979,gutmanmatch1979}.

% The \emph{Hosoya index} $Z(G)$ is defined as the total number of the
 %matchings, i.e.,
 %$$Z(G)=m(G,0)+m(G,1)+\cdots+m(G,\lfloor\frac{n}{2}\rfloor)$$
 The eigenvalues $\lambda_{1}, \lambda_{2},\ldots, \lambda_{n}$ of the adjacency matrix $A(G)$ of $G$ are said to be the eigenvalues of the graph $G$.
 The $energy$ of $G$ is defined as
 \begin{equation}
 E(G)=\sum_{i=1}^{n}|\lambda_{i}|.
 \end{equation}
 The theory of graph energy is well developed nowadays, for details see \cite{gutman2001,GLZ,li&shi&gutman2012}.
 The Coulson integral formula \cite{gutman&polansky1986} plays an important role in the study on graph energy, its version for an acyclic graph $T$ is as follows:
\begin{equation}\label{cformula}
E(T)=\frac{2}{\pi}\int_0^{+\infty}\frac{1}{x^2}\ln\Big[\sum\limits_{t\geq
0} m(T,t)x^{2t}\Big]dx.
\end{equation}
  Motivated by formula \eqref{cformula}, Gutman and Wagner \cite{gutman&Wagner2012} defined the \emph{matching
  energy} of a graph $G$  as
\begin{equation}\label{matchenergyformula}
ME=ME(G)=\frac{2}{\pi}\int^{+\infty}_{0}\frac{1}{x^2}\ln\Big[
\sum\limits_{t\geq 0}m(G,t) x^{2t}\Big]dx. \qquad\qquad
\end{equation}
 Energy and matching energy of graphs are closely related, and they are two quantities of relevance for chemical applications,
 for details see \cite{gutmanmt1975,aihara1976,gutmanmt1977}.

The following result gives an equivalent definition of matching
energy.

\begin{defi} \textnormal{\cite{gutman&Wagner2012}}
 Let $G$ be a graph of order $n$, and let $\mu_1,\mu_2$, $\cdots$, $\mu_n$ be the roots of its matching polynomial. Then
\begin{equation}
ME(G)=\sum\limits_{i=1}^n|\mu_i|.
\end{equation}
\end{defi}

The formula \eqref{matchenergyformula}
% led Gutman and Wagner \cite{gutman&Wagner2012} to introduce
induces a \emph{quasi-order} relation over the set of all graphs on $n$
vertices: if $G_1$ and $G_2$ are two graphs of order $n$, then
\begin{equation}\label{quasi-order} G_1\preceq G_2 \Leftrightarrow m(G_1,t)\leq m(G_2,t) \mbox{ for
all } t=0, 1, \ldots, \lfloor \frac{n}{2}\rfloor.\end{equation}
 If  $G_1\preceq G_2$ and there
exists some $i$ such that $m(G_1,i) < m(G_2, i )$, then we write
$G_1\prec G_2$. Clearly, $$ G_1\prec G_2\Rightarrow
ME(G_1)<ME(G_2).$$

Recall that the \emph{Hosoya index} of a graph $G$ is defined as
$Z(G)=\sum\limits_{t\geq 0} m(G,t)$ \cite{Hosoya}. So we also have
that
 $$ G_1\prec G_2\Rightarrow Z(G_1)<Z(G_2).$$

 The following result gives two fundamental identities for the number of $t$-matchings of a graph \cite{Farrell1979,gutmanmatch1979}.

\begin{lem}\label{twomindentity}
Let $G$ be a graph, $e=uv$ an edge of $G$, and
$N(u)=\{v_1(=v),v_2,\ldots,v_j \}$ the set of all neighbors of $u$
in $G$. Then we have
\begin{equation}\label{removeedge}
m(G,t)=m(G-uv,t)+m(G-u-v,t-1),
\end{equation}
\begin{equation}\label{removevertex}
m(G,t)=m(G-u,t)+\sum_{i=1}^{j}m(G-u-v_i,t-1).
\end{equation}
\end{lem}

From Lemma \ref{twomindentity}, it is easy to get the following
result.
\begin{lem}\textnormal{\cite{gutman&Wagner2012}}\label{subgraph}
Let $G$ be a graph and $e$ one of its edges. Let $G-e$ be the
subgraph obtained from G by deleting the edge $e$. Then $ G-e\prec
G$ and $ME(G-e)<ME(G)$.
\end{lem}

  By Lemma \ref{subgraph}, among all graphs on $n$ vertices, the empty graph $E_n$ without edges and the complete graph $K_n$ have, respectively minimum and maximum matching energy \cite{gutman&Wagner2012}.
  It follows from Eqs. \eqref{cformula} and \eqref{matchenergyformula} that $ME(T)=E(T)$ for any tree $T$ \cite{gutman&Wagner2012}. By using the quasi-order relation, it has  also been obtained
  some results on extremal graphs with respect to matching
  energy among some classes of connected graphs with $n$ vertices. For example, the extremal graphs in connected unicyclic,
  bicyclic graphs were determined by \cite{gutman&Wagner2012} and
  \cite{ji&li&shi2013}, respectively; the minimal graphs among connected
  $k$-cyclic ($k\leq n-4$) graphs and bipartite graphs were
  characterized by \cite{ji&ma2014}; the maximal connected graph with
  given connectivity (resp. chromatic number) was determined by
  \cite{li&yan2013}.

 Let $\mathcal {G}_{n,k}$ be the set of connected graphs of order $n $ $(\geq 2)$ with edge connectivity $k$ ($1\leq k\leq n-1$).
  Let $K_{n-1,1}^k$ be the graph, as shown in Fig. 1, obtained from $K_1\cup K_{n-1}$ by adding $k$ edges between $V(K_1)$ and $V(K_{n-1})$.
  In this paper, we show that $K_{n-1,1}^k$ is the unique graph with maximum matching energy (resp. Hosoya index) in $\mathcal {G}_{n,k}$.

\vspace{10mm}

\setlength{\unitlength}{1mm}
\begin{picture}(40,20)

\put(30,20){\circle{20}}\put(50,20){\circle{20}}

\put(32,22){\circle*{1.0}}\put(32,19){\circle*{1.0}}
\put(32,15){\circle*{1.0}}

\put(48,22){\circle*{1.0}}\put(48,19){\circle*{1.0}}
\put(48,15){\circle*{1.0}}
%%%%vertices%%%
\put(32,22){\line(1,0){16}}
\put(32,15){\line(1,0){16}}\put(32,19){\line(1,0){16}}
\multiput(40,15.5)(0,1.5){3}{\circle*{0.7}}
%%%%line%%%%%

\put(27,20){\makebox(0,0){$K_{n-m}$
}}\put(53,20){\makebox(0,0){$K_m$ }} \put(40,20){\makebox(0,0){$k$
}}
%%%%  graph1 %%%%%

\put(105,20){\circle{20}}

\put(87,20){\circle*{1.5}}

\put(103,20){\circle*{1.0}}\put(103,24){\circle*{1.0}}
\put(103,15){\circle*{1.0}}

%%%%vertices%%%

\put(87,20){\line(1,0){16}}\put(87,20){\line(4,1){16}}\put(87,20){\line(3,-1){16}}
\multiput(100,16.5)(0,1.3){3}{\circle*{0.6}}

%%%%%line%%%
\put(82,20){\makebox(0,0){$K_1$
}}\put(109,20){\makebox(0,0){$K_{n-1}$}}
\put(100,20){\makebox(0,0){$k$ }}
%\put(116,29){\makebox(0,0){$\left.\rule{0mm}{7mm}\right\}$n-9$$}}
%%%%%graph2%%%
\put(40,9){\makebox(0,0){$K^k_{n-m,m}$ $(k\leq m\leq \lfloor
\frac{n}{2}\rfloor)$}} \put(98,9){\makebox(0,0){$K^k_{n-1,1}$ }}
\put(65,0){\makebox(0,0){Fig. 1 \ Graphs $K^k_{n-m,m}$ and
$K^k_{n-1,1}$. }}
%\put(35,7){\makebox(0,0){$S_n^m$ }}\put(85,7){\makebox(0,0){$B_n^m$ }}
\end{picture}
\vspace{5mm}

\section{Main results}

 First we recall some notations. By $\kappa'(G)$ and $\delta(G)$, we denote
 the edge connectivity and the minimum degree of a graph $G$, respectively.
 Let $S$ be a nonempty proper subset of $V$. We use $G[S]$ to denote the subgraph of $G$ induced by $S$. An \emph{edge cut} of $G$, denoted by $\partial(S)$, is a subset of $E(G)$ of the
 form $[S,\bar{S}]$, where $\bar{S}=V\backslash S$. An edge cut $\partial(v)$ ($v\in V$) is
called a \emph{trivial} edge cut. A \emph{$k$-edge cut} is an edge
cut of $k$ elements.
 Let $G\in \mathcal {G}_{n,k}$. Then $G$ must have a $k$-edge cut $\partial(S)$ with $1\leq |S| \leq \lfloor
 \frac{n}{2}\rfloor$.
 %%%such that $G[S]$ and $G[\bar{S}]$ are connected.

\begin{lem}\label{trivial-edge-cut}
  Let $G\ncong  K_{n-1,1}^k$ be a graph in $\mathcal {G}_{n,k}$ with a trivial $k$-edge cut.
  Then $G\prec K_{n-1,1}^k$.
\end{lem}
 \pf Let $\partial(S)$ be a trivial $k$-edge cut of $G$ with $|S|=1$. Since $G\ncong
 K_{n-1,1}^k$, $G[\bar{S}]$ is a  proper subgraph of $K_{n-1}$. Hence $G$ is a proper subgraph of
 $K_{n-1,1}^k$, and so the result follows from Lemma
 \ref{subgraph}.
 \qed

\begin{lem}\label{m&k}
 Let $G \in \mathcal {G}_{n,k}$ be a graph without trivial $k$-edge cuts. Then for any  $k$-edge cut $\partial(S)$ of $G$ with $2\leq |S| \leq \lfloor \frac{n}{2}\rfloor$, we have $|S|\geq k$.
\end{lem}
 \pf For $k\leq 2$, the assertion is trivial, so suppose $k\geq 3$.
 Assume, to the contrary, that $G$ has a $k$-edge cut $\partial(S)$ with $2\leq |S|\leq
 k-1$. By the facts that $\delta(G)\geq \kappa'(G)=k$ and $G$ has no trivial $k$-edge
 cuts, we have $\delta(G)\geq k+1$, and thus $\sum_{v\in
S}d_G(v)\geq |S|(k+1)$. On the other hand, $\sum_{v\in
 S}d_G(v)=2e(G[S])+k\leq |S|(|S|-1)+k$. Therefore, we have
  $|S|(k+1) \leq |S|(|S|-1)+k$, that is, $(|S|-1)(k-|S|)+|S|\leq 0$, which is a contradiction. Therefore the result holds.
\qed

 For $k\leq m\leq \lfloor \frac{n}{2}\rfloor$, let $K_{n-m,m}^k$ be
 the graph, as shown in Fig. 1, obtained from $K_{n-m}\cup K_m$ by adding $k$ independent edges between $V(K_{n-m})$ and $V(K_m)$.
 It is easy to see that $\kappa'(K_{n-m,m}^k)=k$ and $\kappa'(K_{n-1,1}^k)=k$.

 We next show that for a graph $G \in \mathcal {G}_{n,k}$  without trivial $k$-edge cuts, $G\preceq K_{n-m,m}^k$ for some $m$.
Before this,  we introduce a new graph operation as follows.

 let $G_1$ be a graph in $\mathcal {G}_{n,k}$ such that $G_1$ has a $k$-edge
 cut $\partial(S)$ with $G[S]=K_m$, $G[\bar{S}]=K_{n-m}$, and $k\leq m\leq \lfloor \frac{n}{2}\rfloor$. Suppose that $u_1, u_2\in \bar{S}$, $v_1, v_2\in
 S$, $e_1=u_1v_1$, $e_2=u_1v_2$ are two edges of $\partial(S)$, and $u_2$ is not incident with any edge in $\partial(S)$.
 If $G_2$ is obtained from $G_1$ by delating the edge $e_2$ and adding a new edge $e'_2=u_2v_2$, we say that $G_2$ is obtained from $G_1$ by \emph{Operation
 I}, as shown in Fig. 2. Clearly, $G_2 \in \mathcal {G}_{n,k}$.

\vspace{10mm}

\setlength{\unitlength}{1mm}
\begin{picture}(40,20)

\put(30,20){\circle{20}}\put(50,20){\circle{20}}

\put(32,22){\circle*{1.0}}\put(32,19){\circle*{1.0}}\put(32,25){\circle*{1.0}}
\put(32,15){\circle*{1.0}}

\put(48,22){\circle*{1.0}}\put(48,19){\circle*{1.0}}\put(48,25){\circle*{1.0}}
\put(48,15){\circle*{1.0}}
%%%%vertices%%%
\put(32,25){\line(1,0){16}}\put(32,25){\line(5,-1){16}}
\put(32,15){\line(1,0){16}}\put(32,19){\line(1,0){16}}
\multiput(40,15.5)(0,1.5){3}{\circle*{0.7}}
%%%%line%%%%%
\put(40,23){\makebox(0,0){$e_2$}}\put(40,26){\makebox(0,0){$e_1$}}
\put(30,25){\makebox(0,0){$u_1$}}\put(30,22){\makebox(0,0){$u_2$}}
\put(51,25){\makebox(0,0){$v_1$}}\put(51,22){\makebox(0,0){$v_2$}}
\put(27,20){\makebox(0,0){$K_{n-m}$ }}\put(53,19){\makebox(0,0){$K_m$ }}
%%%%  graph1 %%%%%

\put(85,20){\circle{20}}\put(105,20){\circle{20}}

\put(87,22){\circle*{1.0}}\put(87,19){\circle*{1.0}}\put(87,25){\circle*{1.0}}
\put(87,15){\circle*{1.0}}

\put(103,22){\circle*{1.0}}\put(103,19){\circle*{1.0}}\put(103,25){\circle*{1.0}}
\put(103,15){\circle*{1.0}}

%%%%vertices%%%
\put(87,25){\line(1,0){16}}\put(87,22){\line(1,0){16}}
\put(87,15){\line(1,0){16}}\put(87,19){\line(1,0){16}}

\multiput(95,15.5)(0,1.5){3}{\circle*{0.7}}

%%%%%line%%%

\put(93,23){\makebox(0,0){$e'_2$}}\put(93,26){\makebox(0,0){$e_1$}}
\put(83,25){\makebox(0,0){$u_1$}}\put(83,22){\makebox(0,0){$u_2$}}
\put(105,25){\makebox(0,0){$v_1$}}\put(105,22){\makebox(0,0){$v_2$}}

\put(82,20){\makebox(0,0){$K_{n-m}$ }}\put(108,19){\makebox(0,0){$K_m$ }}
\put(65,20){\makebox(0,0){$\longrightarrow$ }}
\put(68,23){\makebox(0,0){\emph{Operation I} }}
\put(45,10){\makebox(0,0){$G_1$ }}
\put(97,10){\makebox(0,0){$G_2$ }}
%\put(116,29){\makebox(0,0){$\left.\rule{0mm}{7mm}\right\}$n-9$$}}
%%%%%graph2%%%
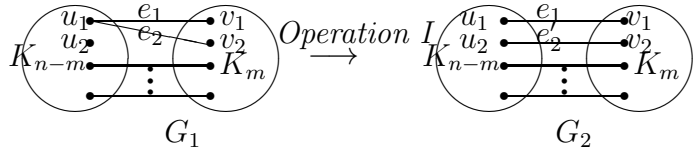
\put(60,5){\makebox(0,0){Fig. 2  The graphs $G_1$ and $G_2$ of
$\mathcal {G}_{n,k}$ in Operation I. }}
%\put(35,7){\makebox(0,0){$S_n^m$ }}\put(85,7){\makebox(0,0){$B_n^m$ }}
\end{picture}

\vspace{-3mm}

\begin{lem}\label{operation}
 If $G_2$ is obtained from $G_1$ by Operation I, then
 $G_1\prec G_2$.
\end{lem}
\pf By formula \eqref{removeedge}, we have
\begin{equation*}
m(G_1,t)=m(G_1-e_2,t)+m(G_1-u_1-v_2,t-1),
\end{equation*}and
\begin{equation*}
m(G_2,t)=m(G_2-e'_2,t)+m(G_2-u_2-v_2,t-1).
\end{equation*}
Note that $G_1-e_2\cong G_2-e'_2$, and $G_1-u_1-v_2$ is isomorphic
to a proper subgraph of $G_2-u_2-v_2$. So, $m(G_1-u_1-v_2,t-1)\leq
m(G_2-u_2-v_2,t-1)$ for all $t$ and
$m(G_1-u_1-v_2,1)<m(G_2-u_2-v_2,1)$. The result thus follows.   \qed

\begin{lem}\label{without trivial-edge-cut}
  Let $G \in \mathcal {G}_{n,k}$ be a graph without trivial $k$-edge cuts.
  Then $G\preceq K_{n-m,m}^k$ for some $m$ with $\max\{k, 2\} \leq m\leq \lfloor \frac{n}{2}\rfloor$.
\end{lem}
 \pf Let $\partial(S)$ be a $k$-edge cut of $G$ with $2\leq |S| \leq \lfloor \frac{n}{2}\rfloor$. Let $|S|=m$. Then $m\geq k$ by Lemma \ref{m&k}.
 Let $G_1$ be the graph obtained from $G$, by adding edges if necessary, such that $G[S]$ and $G[\bar{S}]$
 are complete graphs. Therefore $G\preceq G_1$ by Lemma \ref{subgraph}. If $G_1 \ncong  K_{n-m,m}^k$, then by using
 Operation I repeatedly, we can finally get $K_{n-m,m}^k$ from
 $G_1$. Hence $G_1 \preceq K_{n-m,m}^k$ by Lemma \ref{operation}.
 The proof is thus complete.
  \qed

 In the following, we show that $K_{n-m,m}^k\prec K_{n-1,1}^k$ for $m\geq 2$.

\begin{lem}\label{compare edge}
Suppose $\max\{k, 2\} \leq m\leq \lfloor \frac{n}{2}\rfloor$. Then
$e(K^k_{n-m,m})<e(K^k_{n-1,1}).$
\end{lem}
\pf Note that
\begin{equation*}
\begin{split}
e(K^k_{n-m,m})&=\frac{m(m-1)}{2}+\frac{(n-m)(n-m-1)}{2}+k,
\end{split}
\end{equation*}
and
\begin{equation*} e(K^k_{n-1,1})=\frac{(n-1)(n-2)}{2}+k.
\end{equation*}
 Hence we have
\begin{equation*}
\begin{split}
e(K^k_{n-1,1})-e(K^k_{n-m,m})
&=\frac{n^2-3n+2}{2}-\frac{n^2+2m^2-2mn-n}{2}\\
&=(m-1)(n-m-1)>0.
\end{split}
\end{equation*}
The proof is thus complete. \qed

\begin{lem}\label{n=2m&n=2m+1}
Let $m \geq 1$ be a positive integer. Then we have
\begin{equation}\label{n=2m}
m(K^1_{m,m},t)\leq m(K^1_{2m-1,1},t)   \mbox{ for all } t=0, 1,
\ldots,  m,
\end{equation}
and
\begin{equation}\label{n=2m+1} m(K^1_{m+1,m},t)\leq
m(K^1_{2m,1},t) \mbox{ for all } t=0, 1, \ldots, m.
\end{equation}
\end{lem}

 \pf We apply induction on $m$. For $m=1$ and $m=2$, the assertions are trivial
 since $K^1_{2,2}$ and $K^1_{3,2}$ are proper subgraphs of $K^1_{3,1}$ and $K^1_{4,1}$, respectively.
 So suppose that $m\geq 3$ and Ineqs. \eqref{n=2m} and \eqref{n=2m+1} hold for smaller values of $m$.
  By Lemma \ref{twomindentity}, we obtain that
\begin{equation*}
\begin{split}
m(K^1_{m,m},t)&=m(K^1_{m,m-1},t)+(m-2)m(K^1_{m,m-2},t-1)+m(K_m\cup K_{m-2},t-1)\\
              &=m(K^1_{m,m-1},t)+(m-1)m(K^1_{m,m-2},t-1)-m(K_{m-1}\cup K_{m-3},t-2)\\
              &=m(K^1_{m,m-1},t)-m(K_{m-1}\cup K_{m-3},t-2)+(m-1)[m(K^1_{m-1,m-2},t-1)\\
              &\quad+(m-1)m(K^1_{m-2,m-2},t-2)-m(K_{m-3}\cup K_{m-3},t-3)]\\
              &\leq m(K^1_{m,m-1},t)+(m-1)m(K^1_{m-1,m-2},t-1)+(m-1)^{2}m(K^1_{m-2,m-2},t-2)\\
              &\quad-m(K_{m-1}\cup K_{m-3},t-2),
\end{split}
\end{equation*}and
\begin{equation*}
\begin{split}
m(K^1_{2m-1,1},t)&=m(K^1_{2m-2,1},t)+(2m-3)m(K^1_{2m-3,1},t-1)+m(K_{2m-3},t-1)\\
                 &=m(K^1_{2m-2,1},t)+m(K_{2m-3},t-1)+(2m-3)[m(K^1_{2m-4,1},t-1)\\
                 &\quad+(2m-5)m(K^1_{2m-5,1},t-2)+m(K_{2m-5},t-2)]\\
                 &\geq m(K^1_{2m-2,1},t)+(2m-3)m(K^1_{2m-4,1},t-1)\\
                 &\quad+(2m-3)(2m-5)m(K^1_{2m-5,1},t-2).
\end{split}
\end{equation*}
 By the induction hypothesis, we obtain that
\begin{equation*}
\begin{split}
m(K^1_{m,m-1},t)&\leq m(K^1_{2m-2,1},t),\\
m(K^1_{m-1,m-2},t-1)&\leq m(K^1_{2m-4,1},t-1),\\
m(K^1_{m-2,m-2},t-2)&\leq m(K^1_{2m-5,1},t-2).\\
\end{split}
\end{equation*}
 Since $m\geq 3$, we have that $m-1\leq 2m-3$ and $(m-1)^2\leq(2m-3)(2m-5)$ when
 $n\geq4$. Notice that for $m=3$, $K^1_{m-2,m-2}=K_{m-1}\cup
 K_{m-3}$, and $(m-1)^2-1=(2m-3)(2m-5)$.
 Hence Ineq. (\ref{n=2m}) holds.

By Lemma \ref{twomindentity}, we get that
\begin{equation*}
\begin{split}
m(K^1_{m+1,m},t)&=m(K^1_{m,m},t)+(m-1)m(K^1_{m-1,m},t-1)+m(K_{m-1}\cup K_{m},t-1)\\
              &\leq m(K^1_{m,m},t)+m \cdot m(K^1_{m-1,m},t-1)\\
              &=m(K^1_{m,m},t)+m \cdot [m(K^1_{m-1,m-1},t-1)\\
              &\quad+(m-2)m(K^1_{m-1,m-2},t-2)+m(K_{m-1}\cup K_{m-2},t-2)]\\
              &\leq m(K^1_{m,m},t)+m\cdot [m(K^1_{m-1,m-1},t-1)\\
              &\quad+(m-1)m(K^1_{m-1,m-2},t-2)]\\
              &=m(K^1_{m,m},t)+m \cdot m(K^1_{m-1,m-1},t-1)
              +m(m-1)m(K^1_{m-1,m-2},t-2),
\end{split}
\end{equation*}and
\begin{equation*}
\begin{split}
m(K^1_{2m,1},t)&=m(K^1_{2m-1,1},t)+(2m-2)m(K^1_{2m-2,1},t-1)+m(K_{2m-2},t-1)\\
                 &=m(K^1_{2m-1,1},t)+m(K_{2m-2},t-1)+(2m-2)[m(K^1_{2m-3,1},t-1)\\
                 &\quad+(2m-4)m(K^1_{2m-4,1},t-2)+m(K_{2m-4},t-2)]\\
                 &\geq m(K^1_{2m-1,1},t)+(2m-2)m(K^1_{2m-3,1},t-1)\\
                 &\quad+(2m-2)(2m-4)m(K^1_{2m-4,1},t-2).
\end{split}
\end{equation*}
 By the induction hypothesis and Ineq. (\ref{n=2m}), we have that
\begin{equation*}
\begin{split}
m(K^1_{m,m},t)&\leq m(K^1_{2m-1,1},t)\\
m(K^1_{m-1,m-1},t-1)&\leq m(K^1_{2m-3,1},t-1)\\
m(K^1_{m-1,m-2},t-2)&\leq m(K^1_{2m-4,1},t-2).
\end{split}
\end{equation*}
 Notice that $m\leq 2m-2$ and $m(m-1)\leq (2m-2)(2m-4)$. Therefore Ineq.(\ref{n=2m+1}) holds.

 The proof is thus complete. \qed

\begin{lem}\label{k=1}
 Suppose $2\leq m\leq \lfloor \frac{n}{2}\rfloor$. Then
$$m(K^1_{n-m,m},t)\leq m(K^1_{n-1,1},t) ~\mbox{for all}~ t=0, 1, \ldots, \lfloor \frac{n}{2}\rfloor.$$
\end{lem}

\pf We apply induction on $n$. As the two cases $n=2m$ and $n=2m+1$
 were proved by Lemma \ref{n=2m&n=2m+1}, we proceed to the induction
 step. By Lemma \ref{twomindentity} and the induction hypothesis, we
 have that
\begin{equation*}\label{n-m}
\begin{split}
m(K^1_{n-m,m},t)&=m(K^1_{n-m,m-1},t)+(m-2)m(K^1_{n-m,m-2},t-1)+m(K_{n-m}\cup K_{m-2},t-1)\\
                &\leq m(K^1_{n-m,m-1},t)+(m-1)m(K^1_{n-m,m-2},t-1)\\
                &\leq m(K^1_{n-2,1},t)+(m-1)m(K^1_{n-3,1},t-1)\\
                &=m(K_{n-2},t)+m(K_{n-3},t-1)\\
                &\quad +(m-1)(m(K_{n-3},t-1)+m(K_{n-4},t-2))\\
                &=m(K_{n-2},t)+m\cdot m(K_{n-3},t-1)+(m-1)m(K_{n-4},t-2),
\end{split}
\end{equation*}and
\begin{equation*}
\begin{split}
m(K^1_{n-1,1},t)&=m(K_{n-1},t)+m(K_{n-2},t-1)\\
                &=m(K_{n-2},t)+(n-2)m(K_{n-3},t-1)\\
                &\quad+m(K_{n-3},t-1)+(n-3)m(K_{n-4},t-2)\\
                &=m(K_{n-2},t)+(n-1)m(K_{n-3},t-1)+(n-3)m(K_{n-4},t-2).
\end{split}
\end{equation*}
 Thus the result follows by the fact $m\leq n-2$.\qed

\begin{lem}\label{twopartitionc}
Suppose $k\leq m\leq \lfloor \frac{n}{2}\rfloor$. Then
$$m(K^k_{n-m,m},t)\leq m(K^k_{n-1,1},t) ~\mbox{for all}~ t=0, 1, \ldots, \lfloor \frac{n}{2}\rfloor.$$
\end{lem}

\pf We apply induction on $k$. As the case $k=1$ was proved by
 Lemma \ref{k=1}, we suppose that $k\geq 2$ and the assertion holds for smaller values of $k$.
 By formula \eqref{removeedge}, we have that
 \begin{equation*}
 m(K^k_{n-m,m},t) = m(K^{k-1}_{n-m,m},t)+ m(K^{k-1}_{n-m-1,m-1},t-1),
\end{equation*}
and
\begin{equation*}
m(K^k_{n-1,1},t) = m(K^{k-1}_{n-1,1},t)+ m(K_{n-2},t-1).
\end{equation*}
By the induction hypothesis and Lemma \ref{subgraph}, we obtain that
$m(K^{k-1}_{n-m,m},t)\leq m(K^k_{n-1,1},t)$ and
$m(K^{k-1}_{n-m-1,m-1},t-1) \leq m(K_{n-2},t-1)$. Thus the result
follows. \qed

Together with Lemmas \ref{compare edge} and \ref{twopartitionc}, we
directly obtain the following result.

\begin{cor}\label{compare two graphs}
Suppose $\max\{k, 2\} \leq m\leq \lfloor \frac{n}{2}\rfloor$. Then
$K^k_{n-m,m}\prec K^k_{n-1,1}.$
\end{cor}

\begin{thm}\label{mainresult}
Let $G$ be a graph in $\mathcal {G}_{n,k}$. Then $ME(G)\leq
ME(K^k_{n-1,1})$. The equality holds if and only if $G\cong
K^k_{n-1,1}$.
\end{thm}

\pf  Notice that $K^k_{n-1,1} \in \mathcal {G}_{n,k}$. Let  $G\ncong
K^k_{n-1,1}$ be a graph in $\mathcal {G}_{n,k}$. It suffices to show
that $G\prec K^k_{n-1,1}$. If $G$ has a trivial $k$-edge cut, then
we have $G\prec K_{n-1,1}^k$ by Lemma \ref{trivial-edge-cut}.
Otherwise, by Lemma \ref{without trivial-edge-cut} and Corollary
\ref{compare two graphs}, we obtain that $G\prec K_{n-1,1}^k$ again.
The proof is thus complete. \qed

 By the proof of Theorem \ref{mainresult} and the definition of Hosoya index, we can get the following result on Hosoya index.

 \begin{thm}
Let $G$ be a graph in $\mathcal {G}_{n,k}$. Then $Z(G)\leq
Z(K^k_{n-1,1})$. The equality holds if and only if $G\cong
K^k_{n-1,1}$.
\end{thm}


\begin{thebibliography}{99}

\bibitem{aihara1976}
J. Aihara, A new definition of Dewar-type resonance energies, {\it J. Am. Chem. Soc.}
{\bf 98} (1976) 2750--2758.


\bibitem{graphBondy2008}
J. A. Bondy, U. S. R. Murty, Graph Theory, {\it Springer-Verlag}, Berlin, 2008.



%%%\bibitem{caporossi1999}
%%%G. Caporossi, D. Cvetkovi$\acute{c}$, I. Gutman, P. Hansen, Variable neighborhood search for extremal graphs. 2. Finding graphs with external energy, {\it J. Chem. Inf. Comput. Sci.} {\bf 39} (1999) 984--996.

%%%\bibitem {C.D1980}
%%%D. Cvetkovi\'{c}, M. Doob, H. Sachs, Spectra of Graphs - Theory and
%%%Application, {\it Academic Press}, New York, 1980.


\bibitem{CDGT}
D. Cvetkovi\'{c}, M. Doob, I. Gutman, A. Torga\v{s}ev, Recent
Results in the Theory of Graph Spectra, {\it Elsevier Science
Publishers}, North-Holland, Amsterdam, 1988.

%%%%\bibitem{C.P2010}
%%%%D. Cvetkovi\'{c}, P.Rowlinson, S.Simi\'{c}, An In troduction to the Theory of Graph Spectra, {\it Cabridge University press },Cambridge,2010.

\bibitem{Farrell1979}
E. J. Farrell, An introduction to matching polynomials, {\it J.
Comb. Theory B} {\bf 27} (1979) 75--86.


%%\bibitem {gutman1977}
%%I. Gutman, Acylclic systems with extremal H$\mathrm{\ddot{u}}$ckel $\pi$-electron energy, {\it Theor. Chim. Acta} {\bf 45}(1977), 79--87.

%\bibitem{gutman1982}
%I. Gutman, Correction of the paper ``Graphs with greatest number of
%matchings", {\it Publ. Inst. Math.}(Beograd) {\bf 32} (1982) 61--63.

%\bibitem{gutman1980}
%I. Gutman, Graphs with greatest number of matchings, {\it Publ.
%Inst. Math.} (Beograd) {\bf 27} (1980) 67--76.
\bibitem {gutmanmatch1979}
I. Gutman, The matching polynomial, {\it MATCH Commun. Math. Comput.
Chem.} {\bf 6} (1979) 75--91.

\bibitem {gutman2001}
I. Gutman, \emph{The Energy of a Graph: Old and New Results}, in: A.
Betten, A. Kohnert, R. Laue, A. Wassermann (Eds.), Algebraic
Combinatorics and Applications, {\it Springer-Verlag}, Berlin, 2001,
pp.196--211.

%\bibitem {gutmanCvetkovi1984}
%I. Gutman, D. Cvetkovi\'c, Finding tricyclic graphs with a maximal
%number of matchings  -- another example of computer aided research
%in graph theory, {\it Publications de l'Institut Math\'ematique\/}.
%(Beograd) {\bf 35} (1984) 33--40.


\bibitem{GLZ}
I. Gutman, X. Li, J. Zhang, Graph Energy, in: M. Dehmer, F.
Emmert-Streib (Eds.), {\it Analysis of Complex Networks: From
Biology to Linguistics}, Wiley-VCH, Weinheim, (2009) 145--174.


\bibitem{gutmanmt1975}
I. Gutman, M. Milun, N. Trinajsti$\acute{c}$, Topological definition
of delocalisation energy, {\it MATCH Commun. Math. Comput. Chem.}
{\bf 1} (1975) 171--175.

\bibitem{gutmanmt1977}
I. Gutman, M. Milun, N. Trinajsti$\acute{c}$, Graph theory and molecular orbitals 19, nonparametric resonance energies of arbitrary conjugated systems, {\it J. Am.
Chem. Soc.}{ \bf 99} (1977) 1692--1704.

\bibitem {gutman&polansky1986}
I. Gutman, O.E. Polansky, Mathematical Concepts in Organic
Chemistry, {\it Springer-Verlag}, Berlin, 1986.

\bibitem{gutman&Wagner2012}
I. Gutman and S. Wagner, The matching energy of a graph, {\it
Discrete Appl. Math.} {\bf 160 (15)} (2012) 2177--2187.

\bibitem{Hosoya}
 H. Hosoya, Topological index. A newly proposed quantity characterizing the topological nature of structural isomers of saturated hydrocarbons, {\it Bull. Chem. Soc. Jpn.} {\bf 44} (1971) 2332--2339.

\bibitem{ji&li&shi2013}
S. Ji, X. Li, Y, Shi, The Extremal matching energy of bicyclic
graphs, {\it MATCH Commun. Math. Comput. Chem.} {\bf 70} (2013)
697--706.

\bibitem{ji&ma2014}
S. Ji, H. Ma, The Extremal matching energy of graphs, Ars Combin.,
in press.

\bibitem{li&yan2013}
S. Li, W. Yan, The matching energy of graphs with given parameters,
{\it Discrete Appl. Math.} {\bf 162} (2014) 415--420.
%%\bibitem{koolen&m2001}
%%J.H. Koolen, V. Moulton, Maximal energy graphs, {\it Adv. Appl. Math.} {\bf 26} (2001) 47-52.


\bibitem{li&shi&gutman2012}
X. Li, Y. Shi, I. Gutman, Graph Energy, {\it Springer-Verlag}, New York, 2012.

%%%%%\bibitem {lz2007}
%%%%X. Li, J. Zhang, L. Wang, On bipartite graphs with minimal energy,
%%%%{\it Discrete Appl. Math.} {\bf 157} (2009), 869--873.
\end{thebibliography}
\end{document}